# A note on the sublinear Sobolev inequality

Daniel Spector[*1]

[1]Department of Applied Mathematics, National Chiao Tung University, Hsinchu, Taiwan

May 11, 2016

## Abstract

In this note we show that the classical Sobolev inequality cannot in unrestricted form hold for exponents $p \in (0,1)$.

## 1 Main Result

Let $d \geq 2$. A classical inequality due to Sobolev [10] states that for each $p \in (1,d)$ there exists a universal constant $C = C(p,d) > 0$ such that

$$\|u\|_{L^q(\mathbb{R}^d)} \leq C \|\nabla u\|_{L^p(\mathbb{R}^d;\mathbb{R}^d)} \tag{1.1}$$

for all $u \in L^1_{loc}(\mathbb{R}^d)$ with sufficient decay at infinity and whose distributional derivative $\nabla u \in L^p(\mathbb{R}^d;\mathbb{R}^d)$, where the exponents $p,q$ satisfy the relation

$$\frac{1}{q} = \frac{1}{p} - \frac{1}{d}.$$

The result was subsequently extended to the case $p = 1$ by Gagliardo [3, 4] and Nirenberg [7], who independently recognized the inequality (1.4) as an endpoint of a more general family of interpolation inequalities.

Recently N.C. Phuc [8] has shown that under the further assumptions $\inf_{\mathbb{R}^d} u = 0$ and

$$-div(|Du|^{r-2}Du) \geq 0,$$

for any $r \in (1,2)$, the inequality (1.4) persists for $p \in (r-1,1)$ (see also the work of Buckley and Koskela [1] for their treatment of Sobolev-Poincaré inequalities in the sublinear regime). Here one should be careful to notice both the notion of supersolution and that of gradient require care - as usual for the $r$-Laplace's equation in this regime, supersolutions are based on comparison with solutions, while the gradient is that obtained by the limit of truncations

$$Du(x) := \lim_{k \to \infty} \nabla \min\{u(x), k\}.$$

---

[*]dspector@math.nctu.edu.tw, D.S. supported by MOST 103-2115-M-009-016-MY2

These conditions are readily applicable to a wide class of functions, for example, the $d+3$ parameter family

$$u_{x_0,a,b,c}(x) := \frac{1}{(a+b|x-x_0|)^c},$$

where $x_0 \in \mathbb{R}^d$, $a, b > 0$ and $c > d/p - 1$. This restriction on $c$ ensures $Du \in L^p(\mathbb{R}^d; \mathbb{R}^d)$, while the computation

$$-div(|Du|^{r-2}Du) = cb|cb|^{r-2}\frac{1}{(a+b|x-x_0|)^{(c+1)(r-1)}}\left[\frac{(-c-1)(r-1)b}{(a+b|x-x_0|)} + \frac{d-1}{|x-x_0|}\right]$$

shows that such a function is $r$-superharmonic for any $r$ sufficiently close to one.

**Remark 1.1.** Phuc's results rely on estimates involving the Wolff potential due to Kilpeläinen and Malý [5], and in fact hold for the more general class of $\mathcal{A}$-superharmonic functions (whose prototype is the $r$-Laplacian).

For these classes of functions, Phuc's result represents an improvement over the known results for the sublinear case. One can compare, for instance, with the result of Stein and Weiss [11] whose work implies that for $p \in (1 - 1/d, 1]$ there exists a universal constant $C = C(p, d) > 0$ such that

$$\|u\|_{L^q(\mathbb{R}^d)} \leq C\left(\|\nabla u\|_{L^p(\mathbb{R}^d;\mathbb{R}^d)} + \|(-\Delta)^{1/2}u\|_{L^p(\mathbb{R}^d)}\right), \qquad (1.2)$$

for all $(-\Delta)^{1/2}u \in \mathcal{H}^p(\mathbb{R}^d)$ (which is taken in the sense of an approximation, see [9, p. 123]). Here the fractional Laplacian is defined by its Fourier symbol

$$(-\Delta)^{1/2}u := (2\pi|\xi|\hat{u}(\xi))^{\check{}}.$$

A natural question is whether one has such an improvement to the result on the Hardy spaces without further restriction on the class of functions. In this note, we show that no such improvement is possible in general, which is

**Theorem 1.2.** *There exists a sequence $\{u_n\} \subset W^{1,1}(\mathbb{R}^d)$ with compact support such that for every $q \in (0, \infty)$ and $p \in (0, 1)$*

$$\|u_n\|_{L^q(\mathbb{R}^d)} \to |B(0,1)|^{\frac{1}{q}}$$

*and*

$$\|\nabla u_n\|_{L^p(\mathbb{R}^d;\mathbb{R}^d)} \to 0.$$

*In particular, taking $p, q$ satisfying the scaling relation of Sobolev we see that* (1.4) *fails to hold for general functions with $p \in (0, 1)$.*

**Remark 1.3.** A similar counterexample to the sublinear Sobolev-Poincaré inequality has been exhibited in [1].

**Remark 1.4.** A variation of this argument shows one can take a sequence $\{u_n\} \subset C_c^\infty(\mathbb{R}^d)$ in the preceding theorem.

When $p \in (0, 1 - 1/d)$, it is a result of Krantz [6] (via the atomic decomposition) that the Hardy spaces $\mathcal{H}^p(\mathbb{R}^d)$ continue to provide an appropriate replacement to the inequality (1.4), though the right-hand side of the inequality (1.2) requires further modification. Taking into account Calderón and Zygmund's [2] description of these spaces as higher gradients of harmonic functions, one can deduce the inequality

$$\|u\|_{L^q(\mathbb{R}^d)} \leq C \left( \|\nabla^k u\|_{L^p(\mathbb{R}^d;\mathbb{R}^{kd})} + \|\nabla^{k-1}(-\Delta)^{1/2} u\|_{L^p(\mathbb{R}^d;\mathbb{R}^{(k-1)d})} \right), \quad (1.3)$$

whenever $p \in (\frac{d-1}{k+d-1}, \frac{d-1}{k-1+d-1})$, and where

$$\frac{1}{q} = \frac{1}{p} - \frac{k}{d}.$$

**Remark 1.5.** In principle, Calderón and Zymund had shown that one should require on the right hand side the summation

$$\sum_{i=0}^{k} \|\nabla^{k-i}(-\Delta)^{i/2} u\|_{L^p(\mathbb{R}^d;\mathbb{R}^{(k-i)d})},$$

though we observe that for $i = 2 \ldots k$ one has the pointwise inequality between the Euclidian norms

$$|\nabla^{k-i}(-\Delta)^{i/2} u(x)| \leq C_{p,d} |\nabla^{k-i+2}(-\Delta)^{(i-2)/2} u(x)|$$

since the tensor on the left-hand-side consists of diagonal of the higher order tensor on the right-hand-side.

One can utilize the smooth counterexample referred to in Remark 1.4 to show it is not possible to have a higher order version of Sobolev's inequality whenever $kp < 1$, and so in particular when $d = 2$ this provides a counterexample to a general strengthening of the inequality (1.3) for $p \neq \frac{1}{k}$. Still, it would be of interest to have a counterexample, or family of counterexamples, that works for higher dimensions and all values of $k$ and $p$. More generally, one wonders whether a similar improvement can be made in this regime, whether a higher order version of Sobolev's inequality is true, perhaps restricting the class of functions.

**Open Question 1.6.** Does there exist a universal constant $C = C(p, k, d)$ such that

$$\|u\|_{L^q(\mathbb{R}^d)} \leq C \|\nabla^k u\|_{L^p(\mathbb{R}^d;\mathbb{R}^d)} \quad (1.4)$$

for all $u \in L^1_{loc}(\mathbb{R}^d)$ with sufficient decay at infinity and whose distributional derivative $\nabla^k u \in L^p(\mathbb{R}^d;\mathbb{R}^{kd})$, where the exponents $p, q$ satisfy the relation

$$\frac{1}{q} = \frac{1}{p} - \frac{k}{d},$$

supposing (in an appropriate sense)

$$-div_k(|D^k u|^{r-2} D^k u) \geq 0$$

for some $r > 1$?

*Proof of Theorem 1.2.* We will construct a sequence of piecewise smooth functions $u_n$ converging in pointwise to the characteristic function of the unit ball, while

$$\|\nabla u_n\|_{L^p(\mathbb{R}^d;\mathbb{R}^d)} \to 0,$$

which is a contradiction to the inequality (1.4). Precisely, we consider the sequence

$$u_n(x) = \begin{cases} 1 & \text{if } |x| < 1, \\ n(1 + \frac{1}{n} - |x|) & \text{if } 1 \leq |x| \leq 1 + \frac{1}{n}, \\ 0 & \text{otherwise.} \end{cases}$$

Then, in the sense of distributions we have

$$\nabla u_n(x) = \begin{cases} -n \frac{x}{|x|} & \text{if } 1 < |x| < 1 + \frac{1}{n}, \\ 0 & \text{otherwise,} \end{cases}$$

and so denoting by $\sigma_d$ the measure of the unit sphere in $\mathbb{R}^d$ we find

$$\int_{\mathbb{R}^d} |\nabla u_n|^p \, dx = \sigma_d \int_1^{1+\frac{1}{n}} n^p r^{d-1} \, dr$$
$$= \frac{\sigma_d}{d} n^p \left( (1 + \frac{1}{n})^d - 1 \right).$$

As the function $t \mapsto t^d$ is differentiable at $t = 1$ for any (integer) $d \geq 1$, the right hand side tends to zero as $n \to \infty$. In particular, taking any $p < 1$ and the value of $q$ such that $\frac{1}{q} = \frac{1}{p} - \frac{1}{d}$, we find that (1.4) cannot hold. $\square$

# Acknowledgements

The author is grateful to N.C. Phuc for discussions regarding the problem and Cindy Chen for her assistance in obtaining several of the articles needed in the researching of this paper. The author is supported by the Taiwan Ministry of Science and Technology under research grant MOST 103-2115-M-009-016-MY2 and National Chiao Tung University under research grant 104W986.# References